%Paper: dg-ga/9504003
%From: Xiaowei.Peng@RUBA.RZ.ruhr-uni-bochum.de
%Date: Fri, 28 Apr 1995 13:14:32 +0200

\magnification=\magstep 1
\input amstex
\documentstyle{amsppt}

\def\O{\Omega}
\def\D{\Delta}
\def\A{{\Cal A}}
\def\G{{\Cal G}}
\def\p{\phi}
\def\n{\nabla}
\def\t{\theta}
\def\Ga{\Gamma}
\def\ri{\rightharpoonup}
\def\ZZ{\Bbb Z}
\def\CC{\Bbb C}
\def\NN{\Bbb N}
\def\RR{\Bbb R}
\def\Ima{\hbox{Im}}

\topmatter
\title Variational aspects of the Seiberg-Witten functional
\endtitle
\author J\"urgen Jost, Xiaowei Peng and Guofang Wang
\endauthor
\address Ruhr-Universit\"at Bochum, Fakult\"at f\"ur Mathematik,
44780 Bochum, Germany
\endaddress
\email juergen.jost\@ruba.rz.ruhr-uni-bochum.de
\endemail
\address Ruhr-Universit\"at Bochum, Fakult\"at f\"ur Mathematik,
44780 Bochum, Germany
\endaddress
\email xiaowei.peng\@ruba.rz.ruhr-uni-bochum.de
\endemail
\address Institute of Systems Science, Academia Sinica, 100080 Beijing,
China
\endaddress
\email guofang\@iss06.iss.ac.cn
\endemail
\keywords Seiberg-Witten equation, Seiberg-Witten functional, regularity
and Palais-Smale condition
\endkeywords
\subjclass 58E15, 53C07
\endsubjclass

%\endabstract
\endtopmatter
\document

\head 0. Introduction
\endhead

Recently, Seiberg and Witten (see [SW1], [SW2] and [W]) introduced a new
monopole equation which yields new differential-topological invariants of
four dimensional manifolds, closely related to the Donaldson polynomial
invarints [DK]. This equation has been used to give more elementary proof of
many heorems in gauge theory and to obtain many new results (see [KM], [Le],
[T1], [T2] and [T4]).

These equations are first order equations, but the solutions also satisfy more
general second order equations, in a similar way as (anti) self dual
connections
are solutions of the Yang-Mills equations or holomorphic maps between
K\"ahler manifolds are harmonic. Sometimes, one can use index theorems or
methods from algebraic geometry to construct solutions or to get information
about the space of solutions. Another strategy would be to first solve
the second
order equations and then try to identify conditions under which certain
solutions
of the second order equations actually also solve the first order ones. This
strategy has been successful in many other instances. One of the reasons for
this success usually was the variational structures of the second order
equations, namely, one could identify them as the Euler-Lagrange equations of
some variational integral. The solutions of the first order equations are
absolute minimum of this functional.

We believe that this strategy should also be usefully for the Seiberg-Witten
equations. We therefore study the corresponding variational integral which
we call Seiberg-Witten functional. We show that weak solutions of the
Euler-Lagrange equations are smooth. The main result of our paper is that
the functional satisfies the Palais-Smale condition. As a consquence, one
not only obtains the existence of the minimizers (among which one hopes
to find solutions of the first order equations), but also the existence of
unstable solutions. More precisely, one knows how to deduce from
the Palais-Smale conditon the mountain pass lemma
and a general Ljusternik-Schnirelmann
type theory, and even Morse (-Bott) theory if the functiional is a Morse
(-Bott) functional. We hope that these additional non-minimizing solutions
(which do not solve the first order equations) also carry useful geometric
information and can be used to define more invariants for differential
4-manifolds. This will be studied elsewhere, however.

When we had already nearly completed this paper, the questions that are
answered here were also posed by P. Braam at a conference at the ICTP in
Trieste.

This paper was written while the third auther was a guest of the research
project SFB 237 of the DFG at the Ruhr University Bochum. The first auther
was also supported by the DFG.

\head 1. The Seiberg-Witten equation
\endhead

In order to be able to write down the equation of Seiberg-Witten, we need to
recall the definition of a
$Spin^c$-structure. (For details see [LM]).

Let $(X, g)$ be an oriented, compact Riemannian 4-manifold and $P_{SO(4)}\to
X$ its oriented orthonormal frame bundle.
$Spin^c(4)=Spin(4)\times _{\ZZ_2}U(1)$.
A $Spin^c$-structure is a lift of the
structure group SO(4) to $Spin^c(4)$,
{\it i.e.} there exists a principal $Spin^c$-bundle
$P_{Spin^c(4)}\to X$ such
that there is a bundle map

$$\matrix P_{Spin^c(4)} & \longrightarrow & P_{SO(4)} \cr
&& \cr
\downarrow & & \downarrow \cr
&& \cr
X & \longrightarrow & X \cr
\endmatrix$$

Let $Q=P_{Spin^c(4)}/Spin(4)$ be a principal $U(1)$-bundle.
$ W=P_{Spin^c(4)}\times
_{Spin^c(4)}\CC^4$ and $L=Q\times_{U(1)}\CC$ resp.
is the associated spinor bundle
and the line bundle resp.. $W$ can be decomposed globally as $W^+$ and $W^-$.
Locally,
$$W^{\pm}=S^{\pm}\otimes L^{1/2}.$$
Here $S^{\pm}$ is a spinor bundle w. r. t. a local Spin-struture on $X$. Both
$S^{\pm}$
and $L^{1/1}$ are locally defined.

There exists a Clifford multiplication
$$TX\times W^+\to W^-$$
denoted by $e\cdot \p\in W^-$ for $e\in TX$ and $\p\in W^+$. Here $TX$ is the
tangent
bundle of $X$. A connection on the
bundle $W^+$ can be defined by the Levi-Civita
connection and a connection on $L$.
The ``twisted" Dirac operator $D_A: \Gamma (W^+)
\to\Ga(W^-)$ is defined by
$$D_A=\sum^4_{i=1}e_i\cdot\n_A.$$
Here, $\Ga(W^{\pm})$ is the space of sections of $W^{\pm}$, $\{e_i\}$ is an
orthonormal basis of $TX$ and $\n_A$ is a connection on $W^+$ induced by the
Levi-Civita connection and a connection $A$ on the line bundle $L$.

\definition {Definition 1.1} The Seiberg-Witten equations are
$$
\eqalign{
 D_A\p &=0,\cr
F^+_A &=\frac{i}{4}\langle e_ie_j\p,\p\rangle e^i\wedge e^j,\cr}
\tag 1.1
$$
for $A$ a connection on $L$ and $\p\in\Ga(W^+)$, where $F(A)=-iF_A$ is the
curvature of $A$, $F^+_A$ is the self dual part of $F_A$ and $\{e^i\}$ is the
dual basis of $\{e_i\}$.
\enddefinition

\definition {Definition 1.2} The Seiberg-Witten functional is
$$SW(A,\p)
=\int_X(|D_A|^2+|F^+_A-\frac{i}{4}\langle e_ie_j\p,\p\rangle e^i\wedge
e^j|^2)dvol
\tag 1.2
$$
\enddefinition

The Euler-Lagrange equations of the Seiberg-Witten functional are
$$\eqalign{
D^*_AD_A\p-\frac{i}{2}F^+_A\cdot\p -\frac 18\langle e_ie_j\p,\p\rangle
e_ie_j\p=0,\cr
d^*(F^+_A-\frac{i}{4}\langle e_ie_j\p,\p\rangle e^i\wedge e^j)+
\frac 12\Ima\langle D_A\p,e_i\p\rangle e^i=0.\cr}
\tag 1.3
$$
Here, $D^*_A$ is the formal adjoint operator of $D_A$, $d^*=-*d*$ and $*$
is the Hodge star operator.

It is easy to see that a solution of (1.1) is a solution of (1.3).
In fact, it is clear from the definition of the functional that a solution of
(1.1) is a minimizer of the Seiberg-Witten equation.

The following Weitzenb\"ock formula plays an important role in the
Seiberg-Witten
theory,
$$D^*_AD_A\p=-\Delta _A\p+\frac{s}{4}\p+\frac{i}{2}F_A\cdot\p,
\tag 1.4
$$
where $s$ is the scalar curvature of $(X,g)$.
By this formula, the Seiberg-Witten
functional can be rewritten as
$$SW(A,\p)=\int_X(|\n_A\p|^2+|F^+_A|^2+\frac{s}{4}|\p|^2
+\frac{1}{8}|\p|^4)dvol.
\tag 1.5
$$
And the Euler-Lagrange equation (1.3) can be rewritten as
$$\eqalign{
-\D_A\p+\frac{s}{4}\p+\frac{1}{4}|\p|^2\p=0, \cr
d^*F^+_A+\frac{1}{2}\Ima\langle \n_i\p,\p\rangle e^i=0. \cr}
\tag 1.6
$$
Here $\Delta_A$ is the analyst's Laplacian, the negative Laplacian, and
$\n_i=\n_{e_i}$.

\proclaim{Lemma 1.3} For a smooth solution $(A,\p)$ of equation (1.3)
(or (1.6)),
$|\p|(x)\le\max\{-s_0,0\}$, where $s_0=\min\{s(x)|x\in X\}$.
\endproclaim

\demo{Proof} From the maximum principle.
\qed
\enddemo

\proclaim{Corollary 1.4} If the scalar curvature of $X$ is nonnegative,
then for any smooth
solution $(A,\p)$ of (1.2), $\p\equiv 0$.
\endproclaim

Hence, a solution of (1.2) shares many properties with solutions of the
Seiberg-Witten equation (1.1). In section 3 below, we shall prove Lemma 1.3
for weak solutions of (1.2).

In this paper we shall consider the Seiberg-Witten functional and prove a
compactness theorem. The precise set-up of
the problem will be given in section
2.

Before ending the present section, we want to give the definition of the
{\it Palais-Smale condition}.

Let $M$ be a Banach manifold and $f:M\to \RR$ a smooth functional. Let $G$
be a Lie group acting on $M$ and suppose $f$ is invariant under $G$,
{\it i.e.}
for any $g\in G$ and $x\in X$, $f(gx)=f(x)$. $f$ is said to satisfy
the Palais-Smale
condition: if for any sequence $\{x_i\}_{i\in \NN}$ with
\roster
\item"(i)" $f(x_i)$ is bounded,
\item"(ii)" $df(x_i)\to 0$, as $i\to\infty$,
\endroster
there exists a subsequence (also denoted by $\{x_i\}$)
and a sequence $g_i\in G$ such that $g_ix_i$ converges in $X$ to a critical
point $x$ of $f$, {\it i.e.} $df(x)=0$, with $f(x)=\lim_{i\to\infty}f(x_i)$.

\head 2. The set-up
\endhead

We need to choose a suitable working space to discuss the Seiberg-Witten
functional.

For a vector bundle with a metric over $X$, we define $L^{k,p}(E)$,
the Sobolev
space of sections of $E$ by completing the space $\Ga(E)$ of smooth sections
of $E$ by
$$\|s\|_{L^{k,p}}=(\sum^k_{|\alpha|=0}\int |\n ^{\alpha}s|^p)^{1/p},$$
where $\n $ is a fixed metric connection on $E$. (For details see [P]). For
such Sobolev spaces, the Sobolev embedding theorem and the H\"older
inequality
are valid, {\it i.e.}
$$(\int _X|s|^4)^{1/4}\le c \|s\|_{L^{1.2}}, \hbox{ \qquad for }s\in L^{1,2},
\tag 2.1
$$
and
$$(\int_X|s_1|^2|s_2|^2)^{1/2}
\le (\int _X|s_1|^4)^{1/4}(\int _X|s_2|^4)^{1/4}.
\tag 2.2
$$

Let $L^{1,2}(W^+)$ be the space of sections of the bundle $W^+$ of class
$L^{1,2}$ defined above and ${\Cal  A}^{1,2}=L^{1,2}({\Cal  A})$ the space of
connections of class $L^{1,2}$ defined by a fixed connection $A_0$.
$L^{1,2}(W^+)$ and ${\Cal  A}^{1,2}$ are
Hilbert manifolds.

\proclaim {Lemma 2.1} The Seiberg-Witten functional $SW$ is well defined on
${\Cal  A}^{1,2}\times L^{1,2}(W^+)$ and $SW$ is smooth.
\endproclaim

\demo{Proof} From (2.1) and (2.2),
$$\int_X|A\p|^2\le c \|A\|_{L^{1,2}}\|\p\|_{L^{1,2}}.$$
This implies that $SW$ is well-defined. It is easy to check that $SW$
is smooth.
\enddemo

Now let us to choose a suitable Lie group as a gauge group. First let ${\Cal
G}_0
=\exp(iL^{2,2}(X,\RR))$. We claim that ${\Cal  G}_0$ is a Lie group. Actually,
${\Cal  G}_0$
can be seen as a quotient of $L^{2,2}(X,\RR)$ under an equivalence relation
$\sim$.
$\p_1\sim\p_2$ if and only if $\p_1(x)-\p_2(x)=2\pi n$,
for almost all $x\in X$,
for some integer $n$.
It is clear that $Y=L^{2,2}/\sim$ is a Lie group with the
usual addition of functions. ${\Cal  G}_0$ can be identified with $Y$ by the
exponential map. Hence ${\Cal  G}_0$ is a Lie group with the multiplication of
functions, as the identity component in the continuous case. Then by a
well-known
result about harmonic maps from $X$ into $S^1$ ([EL]), in any component of
$C^{\infty}(X, S^1)$ there exists a unique map $g\in C^{\infty}(X,S^1)$ such
that
$$ d^*(g^{-1}dg)=0,
\tag 2.3
$$
$$ g(x_0)=1, \hbox{ \qquad for  a fixed point \quad}
x_0\in X.
\tag 2.4
$$

\proclaim{Lemma 2.2} Let ${\Cal  G} =\cup g\cdot {\Cal G}_0$. ${\Cal  G}$ is
a Lie group, where the union is over all components of $C^{\infty}(X, S^1)$.
\endproclaim

\demo{Proof} Since ${\Cal  G}_0$ is a Lie group, it suffices to check the
following two points.

(i) For two different components of $C^{\infty}(X,S^1)$ with as $g_1$
and $g_2$
obtained above, $g_1{\Cal  G}_0\cap g_2{\Cal  G}_0=\emptyset$.

Assuming that $g_1{\Cal  G}_0\cap g_2{\Cal  G}_0\ne \emptyset$,
there exist $\varphi_1$
and $\varphi _2\in {\Cal  G}_0$ such that $g_1\varphi_1 =g_2\varphi_2$,
for almost all $x\in X$.
{}From the definition of ${\Cal  G}_0$, there exist $\zeta_1$ and
$\zeta_2\in L^{2,2}(X,\RR)$
such that
$$g_1g_2^{-1}=e^{i(\zeta_2 -\zeta_1)}, $$
for almost all $x\in X$. From (2.3) above, we have
$$d^*d(\zeta_2 -\zeta_1)=0,$$
{\it i.e.} $\zeta_2-\zeta_1 \in L^{2,2}(X,\RR)$ is a harmonic function. Hence
tegether with (2.4), $\zeta_2=\zeta_1$. Therefore $g_1=g_2$, a contradiction.

(ii) The operation of the group is closed.

Let $g_1\varphi_1\in g_1{\Cal  G}_0$,
$g_2\varphi_2 \in g_2{\Cal  G}_0$, then
$$g_1\varphi_1\cdot g_2\varphi_2
=g_1g_2\varphi_1\varphi_2\in g_1g_2{\Cal  G}_0$$
and $g_1g_2$ satisfies (2.3) and (2.4) and is the corresponding $g$ of some
component of $C^{\infty}(X,R)$ as above.
\qed
\enddemo

\proclaim{Lemma 2.3} ${\Cal  G}$ acts smoothly on
$\A ^{1,2}\times L^{1,2}(W^+)$.
\endproclaim

\demo{Proof} $\G$ acts on $\A^{1,2}\times L^{1,2}(W^+)$ as follows
$$g(A,\p)=(g(A),g^{-1}\p),$$
for $(A,\p)\in {\Cal A}^{1,2}\times L^{1,2}(W^+)$ and $g\in\G$, where
$g(A)=A+g^{-1}dg$. It is easy to check that the action is well-defined
and smooth.
\qed
\enddemo

\remark{Remark 2.4} In the non-Abelian case, an element
$g$ of the gauge group
of class $L^{2,2}$ need not act smoothly on $\A^{1,2}$.
\endremark

So on $\A^{1,2}\times L^{1,2}(W^+)$, we can consider the Seiberg-Witten
functional. It is easy to check that (1.2) is equivalent to (1.5) by
an approximation argument. Here, we prefer to use the form (1.5).

\proclaim{Lemma 2.5} The Seiberg-Witten functional $SW$ is coercive, i.e.,
there
exists a constant $c > 0$ such that for each $(A,\p)\in \A^{1,2}\times L^{1,2}
(W^+)$
$$SW(A,\p)\ge c^{-1}(\|g^{-1}\p\|_{L^{1,2}}+\|g(A)\|_{L^{1,2}})-c,$$
for some $g\in{\Cal G}$.
\endproclaim

\demo{Proof} First, we have
$$\eqalign{SW(A,\p)
&= \int(|\n_A\p|^2+|F^+_A|^2+\frac{s}{4}|\p|^2+\frac{1}{8}
|\p|^4) \cr
&=\int\{(|\n_A\p|^2)+\frac{1}{16}|\p|^4 \cr
& +\frac{1}{16}(4s^2+4s|\p|^2+|\p|^4)-\frac{1}{4}s^2\} \cr
&\ge  \int (|\n_A\p|^2+|F^+_A|^2+\frac{1}{16}|\p|^4)-c. \cr}
$$
Using the H\"older inequality, we have
$$(\int|\p|^2)^{1/2}\le c(\int|\p|^4)^{1/4}\le c\int|\p|^4+c.
\tag 2.5
$$
In this paper, $c$ is a constant that may change from term to term.

Also, we have
$$\eqalign{\int_X|\n\p|^2 & \le \int |\n_A\p-(A-A_0)\p|^2 \cr
& \le 2\int(|\n_A\p|^2+|A-A_0|^2|\p|^2) \cr
& \le 2\int|\n_A\p|^2+2(\int |A-A_0|^4)^{1/2}(\int |\p|^4)^{1/2}. \cr}
$$
Together with estimate (2.10) below, it follows that
$$\int|\n\p|^2\le c\int|\n_A \p|^2+c\int|F_A|^2+c.
\leqno(2.6)
$$
Now we esitmate the term containing $g(A)$. There exists a standard method
for the Abelian case. For convenience, we give a complete proof.

Step 1 (gauge fixing). There exists $g_0\in \G_0$ such that
$$d^*(g_0(A)-A_0)=0.
\leqno(2.7)$$
Since $g_0=e^{i\zeta}$ for some $\zeta \in L^{2,2}(X,\RR)$, (2.7) is
equivalent to
$$d^*d\zeta=id^*(A-A_0).
\leqno(2.8)
$$
(2.8) is solvable, for $\int_Xd^*(A-A_0)=0$ and $A\in L^{1,2}$.
By the elliptic estimate, we have
$$
c^{-1}(\|A\|_{L^{1,2}}+\|\p\|_{L^{1,2}})\le
\|g(A)\|_{L^{1,2}}+\|g^{-1}\p\|_{L^{1,2}}
\le c(\|A\|_{L^{1,2}}+\|\p\|_{L^{1,2}}),$$
for some constant $c$.
Hence for simplicity,
we denote $g(A)$ by $A$. So $d^*(A-A_0)=0$.

Step 2 (component fixing). The component group of $C^{\infty}(X,S^1)$ is
isomorphic to $H^1(X, \ZZ)$. For any component of $C^{\infty}(X,S^1)$, there
exists a map $g$ satisfying (2.3) and (2.4). We know
$$g(A)-A_0=g^{-1}dg+A-A_0.$$
The harmonic part of $g(A)-A_0$ is the harmonic part of $A-A_0$ plus the
harmonic part of $g^{-1}dg$. Since the harmonic part of $g^{-1}dg$ belongs to
$H^1(X,\ZZ)$ and the Jacobi torus $H^1(X,\RR)/H^1(X,\ZZ)$ is compact, we can
choose a component such that the harmonic part of $g(A)-A_0$ is bounded.
Since $g$ is harmonic, $d^*(g(A)-A_0)=0$.

Step 3. Using the Hodge decomposition, we have
$$\|g(A)-A_0\|_{L^{1,2}}\le c(\|d^*(g(A)-A_0)\|_{L^2}+\|d(g(A)-A_0)\|_{L^2}+
\|H(g(A)-A_0)\|_{L^2}),$$
where $H(g(A)-A_0)$ is the harmonic part of $g(A)-A_0$. From the preceding
discussion, we obtain
$$\eqalign{\Vert g(A)-A_0\Vert_{L^{1,2}} & \le c\Vert d(g(A)-A_0)
\Vert_{L^2}+c \cr
&\le c\Vert d(g(A))\Vert_{L^2}+c \cr
&\le c\Vert F_A\Vert_{L^2}+c. \cr}
\leqno(2.9)
$$
Therefore,
$$\|g(A)-A_0\|_{L^{1,2}}\le c\|F_A\|_{L^2}+c\le c\|F^+_A\|_{L^2}+c,
\leqno(2.10)
$$
the last inequality is from that $\|F^+_A\|^2_{L^2}-\|F^-_A\|^2_{L^2}$ is
independent of $A$.
Now (2.6) and (2.10) imply
$$SW(A,\p)\ge c^{-1}(\|g^{-1}\p\|_{L^{1,2}}+\|g(A)\|_{L^{1,2}})-c.\quad \qed$$
\enddemo

\head 3. Regularity of weak solutions
\endhead

As in the Yang-Mills case, for the Seiberg-Witten equations there is some kind
of removing singularity theorem (see [Z] and for Yang-Mills see [FU] and [U]).
Actually,
we shall prove in this section that all weak solutions of the second order
equations (1.3) are smooth. This result will be used in the proof of the Main
Theorem.

\proclaim{Theorem 3.1} Let $(A,\p)\in \A^{1,2}\times L^{1,2}(W^+)$ be a weak
solutin of (1.3), i.e. $(A,\p)$ is a critical point of the Seiberg-Witten
functional. Then there exists a gauge transfomation $g\in {\Cal G}$
such that $g(A,\p)=(g(A),g^{-1}\p)$ is smooth.
\endproclaim

First, we show the boundedness of $\| \p\|_{L^{\infty}}$ for a weak solution.

\proclaim{Lemma 3.2}  Let $(A,\p)\in \A^{1,2}\times L^{1,2}(W^+)$
be a weak solution
of (1.3). Then
$$\|\p\|_{L^{\infty}}\le \max \{-\min_{x\in X}s(x),0\}.$$
\endproclaim

\demo{Proof} We use the method of Taubes [T3] to prove the lemma.

Let $s_0=\min_{x\in X}s(x)$.
If $s_0\ge 0$, then from the first equation of (1.6), (recall that (1.3) and
(1.6) are equivalent) we have
$$\int |\n _A\p|^2+\frac{s}{4}|\p|^2+\frac14|\p|^4=0,$$
it follows that $\p=0$. So we may assume $s_0=-1$. Define a test section
$\eta\in L^{1,2}(W^+)$ by

$$ \eta=\cases (|\phi|-1){\phi \over |\phi |}, & \hbox{ for }
|\phi| > 1,\cr 0, & \hbox{ for } |\phi| \le 1. \cr \endcases \tag 3.1$$

Let $\nu=\p/|\p|$ for $ |\p|\ge 1$. It is clear that $|\nu|=1$ and
$$\n\eta =(d|\p|)\nu+(|\p|-1)\n\nu.$$
Since $(A,\p)$ is a weak solution of (1.6), we have
$$\eqalign{
0 & =\int_{\O}\langle \n_A\p,\n_A\eta\rangle +\frac{s}{4}\langle \p,\eta\rangle
+
\frac{1}{4}|\p|^2\langle \p,\eta\rangle  \cr
& =\int_{\O}\langle \n_A\p,\n_A\eta\rangle +\frac 14(|\p|^2+s)(|\p|-1)|\p| \cr
& \ge \int_{\O}\langle \n_A,\n_A\eta\rangle +\frac 14(|\p|^2-1)(|\p|-1)|\p|,
\cr}
$$
for $s\ge -1$. Here, $\O=\{x\in X||\p|(x) > 1\}$. The second term of the
last integration is nonnegative. And the first term is also nonnegative.
In fact, we have, for $x\in\O$,
$$\eqalign{\langle \n_A\p,\n_A\eta\rangle  & =\langle \n_A\p,d|\p|\nu\rangle +
\langle \n_A\p,(|\p|-1)\n_A\nu\rangle
\cr
& =(|\p|-1)|\n_A\nu|^2+\langle \n_A\p,\nu\rangle ^2+(|\p|-1)^2|\n_A\nu|^2 \cr
&+(|\p|-1)\langle \n_A\p,\nu\rangle \langle \nu,\n_A\nu\rangle  \cr
& \ge \frac 12\langle \n_A\p,\nu\rangle ^2+(|\p|-1)|\n_A\nu|^2 \cr
& \ge 0, \cr}
\leqno(3.2)
$$
it follows that the set $\O$ has measure zero, in other words,
$$\|\p\|_{L^{\infty}}\le1.$$
In the general case, {\it i.e.} without the normaliztion $s_0=-1$, the
prceding arguments imply
$$\|\p\|_{L^{\infty}}\le\max\{-s_0,0\}.$$
This completes the proof of the lemma.
\qed
\enddemo

The proof of the Theorem 3.1 is now easy:

We are assuming that $(A,\p)$ is a critical point of $SW$, and thus that
$SW(A,\p)$ is bounded. Lemma 2.5 then implies bounds for
$\|g^{-1}\p\|_{L^{1,2}}$
and $\|g(A)\|_{L^{1,2}}$. Here we also denote $(g(A),g^{-1}\p)$ by $(A,\p)$.

Next, we have
$$ d^*F_A=2d^*F^+_A=-\Ima\langle \n_i^A\p,\p\rangle e^i \hbox{\qquad (from
(1.6))}
\leqno (3.3)
$$
Since $\|\p\|_{L^{\infty}}$ is bounded by Lemma 3.1, this implies
$$ \|d^*F_A\|_{L^2}\le c(\|\n\p\|_{L^2}+\|A\|_{L^2}).
\leqno(3.4)
$$
The second Bianchi identity $dF_A=0$ and the ellipticity of $d+d^*$ imply
$$\Vert F_A\Vert_{L^{1,2}}\le c(\Vert d^*F_A\Vert_{L^2}+\Vert F_A\Vert_{L^2})
\tag 3.5
$$
(3.4), (3.5) and the $L^{1,2}$ estimate for $A$ yield

$$\|A\|_{L^{2.2}}\le c,
\leqno(3.6)
$$
and by the Sobolev embedding theorem then also
$$\|A\|_{L^r}\le c, \hbox{\qquad for  any } r < \infty.
\leqno(3.7)
$$
{}From (1.6), we get , for $1 < p < 2$
$$\|\D\p\|_{L^p}\le c(\|\D_A\p\|_{L^p}
+\|\n A\|_{L^p}+\Vert |A||\n\p|\Vert_{L^p}
+\Vert |A|^2\Vert_{L^p}),$$
and from the H\"older inequality
$$\Vert |a||\n\p|\Vert_{L^p}
\le\Vert A\Vert_{L^{\frac{2r}{2-p}}}\Vert\n\p\Vert_{L^2}.
$$
Thus, by Sobolev's embedding theorem again, $\p\in L^{1,\frac{4p}{4-p}}$
and we may then apply the same kind of argument also for $p=2$ and get
$$\Vert \p\Vert_{L^{2,2}}\le c.$$
The standard bootstrap argument then gives $L^{k,2}$ bounds for $(A,\p)$
for any $k\ge2$, hence smoothness.
\qed

\head 4.The compactness theorem
\endhead

\proclaim{Main Theorem} The Seiberg-Witten functional $SW$ satisfies the
following Palais-Smale condition:

For any sequence $(A_n,\p_n)\in \A^{1,2}\times L^{1.2}(W^+)$ satisfying
\roster
\item"(i)" $dSW(A_n,\p_n)\to 0$ strongly in $\A^{-1,2}\times
L^{-1,2}(W^+)$;
\item"(ii)" $SW(A_n,\p_n)\le c,$ for $n=1,2,\dots,$
\endroster
there exists a subsequence (denoted by $(A_n,\p_n)$) and $g_n\in \G $
such that $g_n(A_n,\p_n)$ converges in $\A^{1,2}\times L^{1,2}(W^+)$ to
a critical point $(A,\p)$ of $SW$ with $SW(A,\p)=\lim_{n\to\infty}SW
(g_n(A,\p))$.
\endproclaim

As we know, the crucial point in the Seiberg-Witten theory and in the
preceding arguments is the boundedness of $\Vert \p\Vert_{L^{\infty}}$.
But for a Palais-Smale sequence (a sequence satisfying (i) and (ii)), there
exists no uniform bound for $\Vert\p_n\Vert_{L^{\infty}}$. This is the main
difficulty we encounter here. Fortunately, we can obtain a weaker bound from
the proof of (3.2) that suffices for showing the Palais-Smale condition.

{\bf Proof of Main Theorem.}

Step 1. By Lemma 2.5, there exist $g\in{\Cal G}$ with
$$\Vert g_n(A_n)\Vert_{L^{1,2}}+\Vert g^{-1}\p\Vert_{L^{1,2}}\le c$$
(independent of n). For simplicity, we denote $g_n(A)$ by $A_n$, and $g^{-1}
\p_n$ by $\p_n$.

{}From Rellich's Theorem and Sobolev's embedding theorem, there exists a
subsequence (also denoted by $(A_n,\p_n)$) such that
\roster
\item"(i)" $A_n\ri A$ weakly in $\A^{1,2}$, and $\p_n\ri\p$ weakly in
$L^{1,2}(W^+)$.
\item"(ii)" $A_n\ri A$ weakly in $\A^{0,4}$, and $\p_n\ri\p$ weakly in
$L^{0,4}(W^+)$.
\item"(iii)" $A_n\to A$ strongly in $\A^{0,p}$, $p < 4$, and $\p_n\to\p$
strongly in $L^{0,p}(W^+)$ $(p < 4)$.
\endroster

Step 2. $(A,\p)$ is a weak solution of (1.6).

For any 1-form $\t\in \A^{1,2}$ (here we abuse the notation a bit),
$$\int\langle F_{A_n}, d\t\rangle +\Ima\langle \n_i^{A_n}\p_n,\p_n\rangle
\langle e^i,\t\rangle =dSW(A_n,\p_n)(\t)=o(1).
\leqno(4.1)
$$
{}From step 1 (i), we know
$$\int \langle F_{A_n},d\t\rangle =\int\langle F_A,d\t\rangle +o(1).
\leqno(4.2)
$$
Now we show that
$$\int \Ima\langle \n_i^{A_n}\p_n,\p_n\rangle \langle e^i,\t\rangle =\int
\Ima\langle \n_i^A\p,\p\rangle \langle e^i,\t\rangle +o(1)
$$
This follows from

$$\eqalign{&\int \Ima\langle \n_i^{A_n}\p_n,\p_n\rangle \langle e^i,\t\rangle
-\int
\Ima\langle \n_i^A\p,\p\rangle \langle e^i,\t\rangle  \cr
&=\int\langle e^i, \t\rangle (\Ima\langle \n_i\p_n+A^i_n\p_n), \p_n\rangle
-\Ima\langle \n_i\p+A^i\p,\p\rangle ) \cr
&=\int\langle e^i, \t\rangle (\Ima\langle \n_i(\p_n-\p)+A^i_n\p_n-A^i\p,
\p_n\rangle -\Ima\langle \n_i\p+A^i\p,\p-\p_n\rangle ) \cr
&=\int\langle e^i, \t\rangle (\Ima\langle \n_i(\p_n-\p)+A^i_n\p-A^i\p,\p\rangle
-\Ima\langle \n_i\p_n+A^i\p,\p\rangle) +o(1) \cr
&\hbox{ \qquad (from step 1 (i) and (iii)) } \cr
&=\int \{\Ima\langle A^i_n(\p_n-\p),\p\rangle +\Ima\langle
(A^i_n-A^i)\p,\p\rangle \}\langle e_i,\t\rangle +o(1) \cr
&=o(1). \hbox{\qquad (from step 1 (iii))} \cr
}
\leqno(4.3)
$$
Here $A-A_0=A^ie^i$ and $A_n-A_0=A^i_ne^i$. From (4.1)--(4.3),
we have for any
$\t$,
$$\int \langle F_A,d\t\rangle +\Ima\langle \n^A_i\p,\p\rangle \langle
e^i,\t\rangle =0,
$$
{\it i.e.} $(A,\p)$ satisfies weakly the second equation of (1.6). Using the
same argument, we can show that $(A,\p)$ satisfies weakly
the first equation of
(1.6). That is, $(A,\p)$ is a weak solution of (1.6).

Hence from Theorem 3.1, there exists $g\in {\Cal G}$ such that $g(A,\p)$ is a
smooth
solution of (1.6) and $|\p(x)|\le s_0$
(recall that $s_0=\max\{-\min_{x\in X}s(x),0\}$).
So we may assume that $(A,\p)$ is a smooth solution.

Step 3. As in the proof of Lemma 3.2, we assume $s_0:=\min_{x\in X}s(x)=-1$.
Set
$$\O_n:=\{x\in X||\p_n| > 1\},
$$
and $\nu_n=\p_n/|\p_n|$ for $x\in \O_n$.

\remark {Claim} $$\int_{\O_n}|\langle \n_{A_n}\p_n,\nu_n\rangle |^2\to 0,
\hbox{ as } n\to\infty.
\leqno(4.4)$$
\endremark

As in the proof of lemma 3.2, we have (see (3.2))
$$\int_{\O_n}|\langle \n_{A_n}\p_n, \nu_n\rangle |^2\le 2
dSW(A_n,\p_n)(\eta_n)\le 2 \Vert
dSW(A_n,\p_n)\Vert_{L^{-1,2}}\Vert\eta_n\Vert_{L^{1,2}}.
$$
Here $\eta_n$ is defined as $\eta$ in the proof of Lemma 3.2, namely,

$$ \eta_n=\cases (|\phi_n|-1){\phi_n \over |\phi_n |}, & \hbox{ for }
|\phi_n| > 1, \cr
0, & \hbox{ for } |\phi_n| \le 1. \cr \endcases$$

It is sufficient to show that $\Vert\eta_n\Vert_{L^{1,2}}\le c$ for a constant
$c$.
{}From the definition of $\eta_n$, we have
$$\Vert\eta_n\Vert^2_{L^2}=\int_{\O_n}(|\p_n|-1)^2\le\int_X|\p_n|^2+c\le
\Vert\p_n\Vert^2_{L^{1,2}}+c\le c,
$$
and
$$\eqalign{
\Vert\n\eta_n\Vert^2_{L^2}
&= \Vert d|\p_n|\nu_n +(|\p_n|-1)\n\nu_n\Vert^2_{L^2(\O_n)} \cr
& \le
2\Vert\n\p_n\Vert^2_{L^2}+8\Vert\frac{|\p|-1}{|\p|}\Vert_{L^{\infty}(\O_n)}
\Vert\n\p_n\Vert^2_{L^2} \cr
& \le c\Vert\n\p_n\Vert^2_{L^2}\le c. \cr}
$$

Step 4. $A_n\to A$ strongly in $\A^{1,2}$.

{}From step 1 (i) and (4.1), we have
$$\eqalign{
\Vert F_{A_n}-F_A\Vert^2_{L^2} &= \int _X \langle d(A_n-A),d(A_n-A)\rangle  \cr
& = \int_X \langle dA_n,d(A_n-A)+\int_X\langle dA,d(A_n-A)\rangle  \cr
& =-\int_X\Ima\langle \n^{A_n}_i\p_n,\p_n\rangle \langle e^i, A_n-A\rangle
+o(1) \cr}
\leqno(4.5)
$$

To use step 3, we decompose $X$ as $\O_n$ and $X\backslash\O_n$.
On $X\backslash\O_n$,
$|\p_n|\le 1$.
So it is clear that
$$\eqalign{
& -\int_{X\backslash\O_n}\Ima\langle \n^{A_n}_i\p_n,\p_n\rangle \langle e^i,
A_n-A\rangle  \cr
\le & \int_{X\backslash\O_n}(|\n\p_n|+|A_n|)|A_n-A|^2 \cr
\le & (\int_X (|\n\p_n|^2+|A_n|^2 )^{1/2}\int_X|A_n-A|^2 \cr}
\leqno(4.6)
$$
On the other hand,
$$\eqalign{
& -\int_{\O_n}\Ima\langle \n^{A_n}_i\p_n,\p_n\rangle \langle e^i, A_n-A\rangle
\cr
\le & \int_{\O_n}|\langle \n_i^{A_n}\p_n,\nu_n\rangle ||\p_n||A_n-A| \cr
\le & (-\int_{\O_n}|\langle \n_i^{A_n}, \nu_n\rangle
|^2)^{1/2}\Vert\p_n\Vert_{L^4}
(\Vert A_n\Vert_{L^4}+\Vert A\Vert_{L4}) \cr
& \to 0, \hbox{\qquad as }n\to\infty \hbox{\qquad (using step 3).} \cr }
\leqno(4.7)
$$
(4.5), (4.6) and (4.7) imply
$$\Vert A_n-A\Vert_{L^{1,2}}\to 0,\hbox{\qquad as \quad} n\to\infty.$$

Step 5. $\p_n\to\p$ strongly in $L^{1,2}(W^+)$.

First, from step 1 (i), we have
$$\eqalign{
\Vert\n\p_n-\n\p\Vert^2_{L^2} & = \int_X(\langle \n\p_n,\n(\p_n-\p)\rangle +
\langle \n\p,\n(\p_n-\p)\rangle ) \cr
& = \int _X\langle \n\p_n,\n(\p_n-\p)\rangle +o(1). \cr}
\leqno(4.8)
$$
$\n\p_n=\n_{A_n}\p_n-A_n\p_n$, so we have
$$\eqalign{
&\int_X\langle \n\p_n,\n(\p_n-\p)\rangle  \cr
=& \int_X\langle \n_{A_n}\p_n, \n_{A_n}(\p_n-\p)\rangle
-\int_X\langle A_n\p_n,\n(\p_n-\p)\rangle  \cr
& -\int_X\langle \n\p_n, A_n(\p_n-\p)\rangle
+\int_X\langle A_n\p_n,A_n(\p_n-\p)\rangle . \cr}
\leqno(4.9)
$$
{}From (i) of the Palais-Smale condition,
we know
$$\eqalign{
& \int_X\langle \n_{A_n}\p_n, \n_{A_n}(\p_n-\p)\rangle  \cr
= & -\int _X\frac 14(|\p_n|^2+s)\langle \p_n,\p_n-\p\rangle
+dSW(A_n,\p_n)(\p_n-\p)
\cr
= & -\int _X\frac
14(|\p_n|^2+s_0)\langle \p_n,\p_n-\p\rangle +\int_X(s_0-s)\langle
\p_n,\p_n-\p\rangle +o(1) \cr
= & -\int _X\frac
14(|\p_n|^2+s_0)|\p_n-\p|^2-\int_X\frac
14(|\p_n|^2+s_0)\langle \p,\p_n-\p\rangle +o(1)  \cr
& \hbox{\quad (since } \p_n \to\p \hbox{ strongly in }
L^2) \cr
\le & -\frac
14(|\p_n|^2+s_0)\langle \p,\p_n-\p\rangle +o(1)  \cr
&\hbox{\qquad (since the first term is nonpositive.)} \cr
\le &\frac 14 (\int_X||\p_n|^2+s_0|^2)^{1/2}\int_X|\p_n-\p|^2+o(1) \cr
& \hbox{\qquad (since } |\p|\le s_0 \hbox{ by  Lemma 3.1}) \cr
= & o(1), \cr}
\leqno(4.10)
$$
as $n\to\infty$.

The other three terms can be estimated using step 4.
For example, the second term
$$\eqalign{
& \int _X\langle A_n\p_n,\n(\p_n-\p)\rangle  \cr
= & \int _X\langle (A_n-A)\p_n,\n(\p_n-\p)\rangle +\langle
A\p_n,\n(\p_n-\p)\rangle  \cr
\le & \Vert A_n-A\Vert_{L^4}\|\p_n\Vert_{L^4}\Vert\p_n-\p\Vert_{L^{1,2}}
\cr
&+\int_X \langle A(\p_n-\p),\n(\p_n-\p)\rangle +\langle A\p,\n(\p_n-\p)\rangle
\cr
\le & \Vert A_n-A\Vert_{L^4}\|\p_n\Vert_{L^4}\Vert\p_n-\p\Vert_{L^{1,2}}
\cr
&+\Vert\p_n-\p\Vert_{L^2}\Vert\n(\p_n-\p)\Vert_{L^2}+
\Vert A\p\Vert_{L^{1,2}}\Vert\p_n-\p\Vert_{L^2} \cr
= & o(1). \cr
}
$$

Therefore, we show that $\Vert\p_n-\p\Vert_{L^{1,2}}\to 0$ as $ n\to\infty$.
Together with step 4, this completes the proof of the Main Theorem.


\Refs

\refstyle{A}
\widestnumber\key{SW2}

\ref\key DK
\by S. K. Donaldson and P. B. Kronheimer
\book The geometry of four-manifolds
\bookinfo Oxford Science Publications
\yr 1990
\endref


\ref\key EL
\by J. Eells and L. Lemaire
\paper Another report on harmonic maps
\jour Bull. London Math. Soc.
\vol 20
\yr 1988
\pages 385--524
\endref

\ref\key FU
\by D. Freed and K. K. Uhlenbeck
\book Instantons and four manifolds
\bookinfo Springer-Verlag New York
\yr 1984
\endref

\ref\key KM
\by P. B. Kronheimer and T. S. Mrowka
\paper The genus of embedded surfaces in the projective space
\jour Math. Res. Letters
\vol 1
\yr 1994
\pages 797--808
\endref

\ref\key Le
\by C. LeBrun
\paper Einstein metrics and Mostow rigidity
\jour  Math. Res. Letters
\vol 2
\yr 1995
\pages 1--8
\endref

\ref\key LM
\by  H. B. Lawson and M.-L. Michelsohn
\book  Spin geometry
\bookinfo Priceton, New Jersey
\yr 1989
\endref

\ref\key P
\by R. S. Palais
\book  Foundations of global non-linear analysis
\bookinfo W. A. Benjamin, Inc. New York
\yr 1968
\endref

\ref\key SW1
\by N. Seiberg and E. Witten
\paper Electromagnetic duality, monopole condensation and confinement in
N=2 supersymmetric Yang-Mills theory
\jour Nucl. Phys.
\vol B426
\yr 1994
\pages 19--52
\endref

\ref\key SW2
\by N. Seiberg and E. Witten
\paper Monopoles, duality and chiral symmetry breaking in N=2
supersymmetry QCD
\jour Nucl. Phys.
\vol B431
\yr 1994
\pages 581--640
\endref

\ref\key T1
\by C. H. Taubes
\paper The Seiberg-Witten invariants and symplectic forms
\jour Math. Res. letters
\vol 1
\yr 1994
\pages 809--822
\endref

\ref\key T2
\bysame
\paper The Seiberg-Witten invariants and the Gromov invariants
\jour
\yr
\pages
\endref

\ref\key T3
\bysame
\paper On the equivalence of the first and second order equations
for gauge theories
\jour Comm. Math. Phys.
\vol 75
\yr 1980
\pages 207--227
\endref

\ref\key T4
\bysame
\paper More constraints on symplectic manifolds from Seiberg-Witten
invariants
\jour  Math. Res. Letters
\vol 2
\yr 1995
\pages 9--14
\endref

\ref\key U
\by K. K. Uhlenbeck
\paper Removable singularities in Yang-Mills fields
\jour Comm. Math. Phys
\vol 83
\yr 1982
\pages 31--42
\endref

\ref\key W
\by E. Witten
\paper Monopoles and 4-manifolds
\jour Math. Res. Letters
\vol 1
\yr 1994
\pages 769--796
\endref

\ref\key Z
\by J. Zhou
\paper A vanishing theorem for Seiberg-Witten invariants
\jour
\yr
\pages
\endref

\endRefs
\enddocument